\documentclass{article}
\usepackage{amsmath,mathrsfs,amssymb,amsthm,amscd}
\usepackage[]{fontenc}
\usepackage[all]{xy}
\usepackage{hyperref}
\usepackage{graphics}
\usepackage{a4wide}
\usepackage[dvips]{graphicx,color}

\newcounter{EQNR}
\renewcommand{\contentsname}{Contents\\
{\footnotesize\normalfont(A table of contents should normally not
be included)}}

\begin{document}

\title{On the appearance of Eisenstein series through degeneration}
\author{}
\date{D. Garbin, J. Jorgenson and M. Munn}

\maketitle

\begin{abstract}\noindent
Let $\Gamma$ be a Fuchsian group of the first kind acting on the
hyperbolic upper half plane $\mathbb H$, and let $M = \Gamma
\backslash \mathbb H$ be the associated finite volume hyperbolic
Riemann surface.  If $\gamma$ is parabolic, there is an associated
(parabolic) Eisenstein series, which, by now, is a classical part
of mathematical literature. If $\gamma$ is hyperbolic, then,
following ideas due to Kudla-Millson, there is a corresponding
hyperbolic Eisenstein series.  In this article, we study the
limiting behavior of parabolic and hyperbolic Eisenstein series on
a degenerating family of finite volume hyperbolic Riemann
surfaces. In particular, we prove the following result. If $\gamma
\in \Gamma$ corresponds to a degenerating hyperbolic element, then
a multiple of the associated hyperbolic Eisenstein series
converges to parabolic Eisenstein series on the limit surface.
\end{abstract}

\section{Introduction}

\begin{nn}\label{1.1}
\textbf{Spectral expansions.} Let $M = \Gamma \backslash \mathbb
H$ be a finite volume hyperbolic Riemann surface, realized as the
quotient of the hyperbolic upper half plane $\mathbb  H$ by a
discrete subgroup $\Gamma$ of $\textrm{PSL}_{2}(\mathbb  R)$.  Let
$\Delta_{M}$ denote the Laplacian, associated to the hyperbolic
metric, which acts on the space of smooth functions on $M$.  For
the sake of our discussion, consider the corresponding heat kernel
$K_{M}(t;z,w)$, which is a function of $t \in \mathbb R^{+}$ and
$z, w \in M$.  If $M$ is compact, then the heat kernel admits the
spectral expansion
\begin{equation}\label{heatcompact}
K_{M}(t;z,w)=\sum\limits_{n=0}^{\infty}e^{-\lambda_{M,n}t}\phi_{M,n}(z)
\overline{\phi_{M,n}(w)}
\end{equation}
where $\{\phi_{M,n}\}$ is a complete orthonormal basis of
eigenfunctions of the $\Delta_{M}$ with corresponding
(non-negative) eigenvalues $\lambda_{M,n}$.  If $M$ is
non-compact, then the spectral expansion of the heat kernel takes
a very different form.  More specifically, let $\{P\}$ denote the
finite set of $\Gamma$-inequivalent cusps, and $E_{{\rm
par};M,P}(s,z)$ be the (parabolic) Eisenstein series on $M$ corresponding
to $P$.  Then, the spectral expansion of the heat kernel on $M$ is
the identity
\begin{equation}\label{heatnoncompact}
\begin{array}{l}\displaystyle
K_{M}(t;z,w)=\sum\limits_{n=0}^{\infty}e^{-\lambda_{M,n}t}\phi_{M,n}(z)
\overline{\phi_{M,n}(w)} \\[3mm]\displaystyle \hskip 20mm +
\,\,\,\,
\frac{1}{4\pi}\sum\limits_{P}\int\limits_{-\infty}^{\infty}e^{-(r^{2}+1/4)t}
E_{{\rm par};M,P}(1/2+ir,z) \overline{E_{{\rm
par};M,P}(1/2+ir,w)}dr\,.
\end{array}
\end{equation}

\vskip .10in Recall now that any finite volume hyperbolic Riemann
surface $M_{0}$ can be realized as one component of a degenerating
sequence of compact hyperbolic Riemann surfaces $M_{\ell}$. In
this setting, it has been shown that the hyperbolic heat kernels
on $M_{\ell}$ converge to the hyperbolic heat kernel on $M_{0}$;
see \cite{JLu1} and \cite{JoLu}. With the heat kernel convergence result in mind,
one immediately has from (\ref{heatcompact}) and
(\ref{heatnoncompact}) the following natural question:  How does
one see the emergence of the Eisenstein series in
(\ref{heatnoncompact}) through degeneration? More precisely, does
there exist a naturally defined sequence of functions
$h_{\ell}(s,z)$ on $M_{\ell}$ which converges to the Eisenstein
series $E_{{\rm par};M_{0},P}(s,z)$ on $M_{0}$?
\end{nn}

\begin{nn}\label{1.2}
\textbf{Spectral theory on degenerating Riemann surfaces.} The
problem of studying asymptotic behavior of spectral theory on
degenerating Riemann surfaces of finite volume has received
considerable attention in the literature.  In \cite{HejhalAMS},
Hejhal developed the theory of degenerating $b$-groups and
obtained, among other results, the lead asymptotics of spectral
counting functions.  An improvement of the error term in the case
the degenerating surfaces are compact was proved in \cite{JiZ}.
From \cite{Fay}, one has a construction of degenerating hyperbolic
Riemann surfaces of finite volume by first constructing families
of degenerating algebraic curves, from which one can utilize the
uniformization theorem in order to obtain degenerating families of
Riemann surfaces of finite volume.  In \cite{Jorgen}, the approach
from \cite{Fay} is used to study spectral invariants associated to
the canonical and Arakelov metrics.  Beginning in \cite{JLu1},
Huntley, Jorgenson and Lundelius used the methodology from
\cite{Fay} to study hyperbolic spectral theory through
degeneration.  These authors obtained numerous result, including:
Convergence of heat kernels \cite{JLu1}; asymptotic behavior of
heat traces and Selberg zeta functions \cite{JoLu2}; convergence
of relative spectral functions \cite{JLu2}; asymptotic behavior of
counting functions \cite{JoLu}; asymptotic behavior of weighted
counting functions (Riesz sums) \cite{HJL}.  In all these
articles, the results apply to non-compact degenerating families
as well as compact families. Further results concerning eigenvalue
and eigenfunction convergence have been obtained by Judge in
\cite{JuI} and \cite{JuII}, and Wolpert used degenerating
techniques to study the problem of existence of $L^2$
eigenfunctions on general finite volume hyperbolic Riemann
surfaces. More recently, in \cite{JK10}, it was shown that one can
use the results from \cite{JoLu2} to prove results for other
metrics, namely it was shown that the metric on Teichm\"uller
space induced from the canonical metric is not complete.

\vskip .10in In brief, there is a vast literature addressing
problems in the study of spectral theory on degenerating finite
volume Riemann surfaces.  Further problems exist, and as
mathematical development demonstrates, new results are interesting
for their own sake as well as for potential applications to other
fields.

\end{nn}

\begin{nn}\label{1.3}
\textbf{The main results.}
Throughout this article we refer to the Eisenstein series
$E_{{\rm par};M,P}(s,z)$ in (\ref{heatnoncompact}) as parabolic Eisenstein series since each
such series is associated to a parabolic element of the uniformizing
group $\Gamma$.  In \cite{KM}, the authors defined a hyperbolic Eisenstein
series $E_{{\rm hyp};M,\gamma}(s,z)$ associated to any hyperbolic element
$\gamma \in \Gamma$.  We will summarize both definitions in sections 2.3 and 2.4.
In addition, as we will recall below, a degenerating family of hyperbolic
Riemann surfaces $M_{\ell}$ has two types of hyperbolic elements: Non-degenerating
elements, which are those that converge to hyperbolic elements in the Fuchsian group
of the limit
surface, and degenerating elements, which are those whose associated geodesics have
lengths that converge to zero.

\vskip .10in  Precise definitions and references to all concepts will be given
in section 2 below.  However, with these comments made, we are able to state
the main result of the paper.

\vskip .10in {\bf Main Theorem.} {\it Let $M_{\ell}$ be a
degenerating family of hyperbolic Riemann surfaces of finite
volume, with limit surface $M_{0}$.}
\begin{list}{}
\item{i)} {\it Let $E_{{\rm hyp};M_{\ell},\gamma}(s,z)$ be the
hyperbolic Eisenstein series on $M_{\ell}$ associated to the
hyperbolic element $\gamma$. If $\gamma$ corresponds to a
non-degenerating hyperbolic element, then}
$$
\lim\limits_{\ell_{\gamma} \rightarrow 0} E_{{\rm
hyp};M_{\ell},\gamma}(s,z) = E_{{\rm hyp};M_{0},\gamma}(s,z)\,.
$$
\item{ii)} {\it Let $E_{{\rm par};M_{\ell},P}(s,z)$ be the parabolic
Eisenstein series on $M_{\ell}$ associated to the cusp $P$.  Then}
$$
\lim\limits_{\ell_{\gamma} \rightarrow 0} E_{{\rm
par};M_{\ell},P}(s,z) = E_{{\rm par};M_{0},P}(s,z)\,.
$$
\item{iii)} {\it Let $E_{{\rm hyp};M_{\ell},\gamma}(s,z)$ be the
hyperbolic Eisenstein series on $M_{\ell}$ associated to the
hyperbolic element $\gamma$, whose geodesic has length $\ell_{\gamma}$.
If $\gamma$ corresponds to a degenerating hyperbolic
element which results in the new cusp $P$, then}
$$
\lim\limits_{\ell_{\gamma} \rightarrow 0}
\ell_{\gamma}^{-s}E_{{\rm hyp};M_{\ell},\gamma}(s,z) = E_{{\rm
par};M_{0},P}(s,z)\,.
$$
\end{list}
{\it In all instances, the convergence is uniform on compact
subsets of $M_{0}$ bounded away from the developing cusps, and in
half-planes of the form ${\rm Re}(s) \geq 1 + \delta$ for any $\delta > 0$.}
\end{nn}

\vskip .10in
Part (iii) answers the question posed above, namely to determine a naturally
defined sequence of functions on a degenerating family $M_{\ell}$ of hyperbolic Riemann
surfaces whose limit is the parabolic Eisenstein series associated to the
newly developed cusps.

\vskip .10in {\bf Explanation of (iii).} In order to keep the
statement of part (iii) manageable, we employed a slight abuse of
notation, which we now explain.  If $M_{\ell}$ has a single
pinching geodesic which is separating, then the limit surface
$M_{0}$ has two components, which we denote by $M_{0;1}$ and
$M_{0;2}$, each with a single newly formed cusp, denoted by
$P_{1}$ and $P_{2}$. In this case, the right-hand-side of (iii)
depends on the location of the point $z$:  If $z \in M_{0;1}$
(resp. $z \in M_{0;2}$), then the right-hand-side of (iii)
signifies the function $E_{{\rm par};M_{0;1},P_{1}}(s,z)$ (resp.
$E_{{\rm par};M_{0;2};P_{2}}(s,z)$).  If $M_{\ell}$ has a single
pinching geodesic which is non-separating, then the limit surface
$M_{0}$ has one component with two newly formed cusps, denoted by
$P_{1}$ and $P_{2}$.  In this case, the right-hand-side of (iii)
signifies the function $E_{{\rm par};M_{0},P_{1}}(s,z) + E_{{\rm
par};M_{0},P_{2}}(s,z)$.  To consider the general case when
$M_{\ell}$ has several pinching geodesics, then one simply
iterates the results from the Main Theorem by pinching one
geodesic at a time.

\begin{nn}\label{1.4}
\textbf{Outline of the paper.} In section~2, we establish notation
and recall various known results. Most important for our results
are the definitions of parabolic and hyperbolic counting
functions, and the realization that parabolic and hyperbolic
Eisenstein series can be expressed as Stieltjes integrals of these
counting functions. In section~3, we study the asymptotic behavior
of the counting functions from section~2 through degeneration.
With these results, we conclude by proving the Main Theorem in
section~4.
\end{nn}

\section{Background material}

\begin{nn}\label{2.1}
\textbf{Basic notation.} Let $M$ be a finite volume hyperbolic
Riemann surface.  By this we mean there exists a Fuchsian group of
the first kind $\Gamma$ acting on the hyperbolic upper half plane
$\mathbb H$ such that $M$ and $\Gamma \backslash \mathbb H$ are
isometric.  Hence, we write $M = \Gamma \backslash \mathbb H$. As
is common, we realize $\mathbb H$ as $\{z \in \mathbb C \,\,|\,\,
\textrm{\rm Im}(z) > 0\}$. Writing $z = x+iy$, then the hyperbolic
metric $\mu_{\textrm{hyp}}$ and hyperbolic Laplacian
$\Delta_{\textrm{hyp}}$ can be expressed as
$$
\mu_{\textrm{hyp}} = \frac{dx^2+dy^2}{y^2}
\,\,\,\,\,\textrm{\rm and}\,\,\,\,\,
\Delta_{\textrm{hyp}} =- y^2 \left( \frac{\partial^2}{\partial x^2}
+ \frac{\partial^2}{\partial y^2} \right).
$$
Under the change of coordinates
$x=e^{\rho}\cos\theta$ and $y=e^{\rho}\sin \theta$, the hyperbolic metric
and hyperbolic Laplacian are given by
$$
\mu_{\textrm{hyp}} =
\frac{d\rho^2+ d\theta^2}{\sin^2 \theta} \,\,\,\,\,\textrm{\rm and}\,\,\,\,\,
\Delta_{\textrm{hyp}} =- \sin^2 \theta \left( \frac{\partial^2}{\partial \rho^2}
+ \frac{\partial^2}{\partial \theta^2} \right).
$$
In a slight abuse of notation, we will at times in this
article identify $M$ with a fundamental domain (say, a Ford
domain, bounded by geodesic paths) and identify points on $M$ with
their pre-images in $\mathbb{H}$.
\end{nn}

\begin{nn}\label{2.2}
\textbf{Counting functions.}
\begin{figure}
\begin{center}
\includegraphics[width=0.6\textwidth]{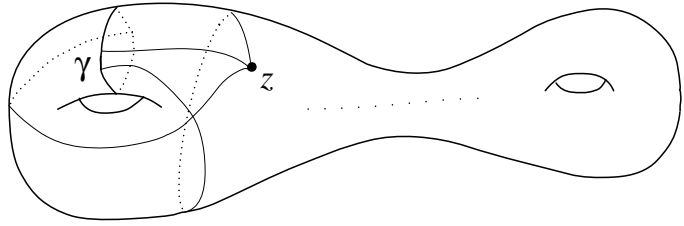}
\caption{Geodesic paths from a point to a closed geodesic}
\end{center}
\end{figure}
Let $\gamma \in \Gamma$ be a primite hyperbolic element. As usual,
primitive means that if $\gamma_{0}\in\Gamma$ and $\gamma_{0}^{n}
= \gamma$ for some integer $n$, then $n = \pm 1$.  By hyperbolic,
one means that $\gamma$ can be conjugated in $\textrm{\rm
PSL}_{2}(\mathbb R)$ to a non-identity diagonal matrix, which we
write as
$$
\gamma = \left(\begin{array}{ll} e^{\ell_{\gamma}/2}& 0 \\ 0 & e^{-\ell_{\gamma}/2}\end{array}
\right),
$$
where $\ell_{\gamma}$ denotes the length of the closed geodesic on $M$ in the homotopy
class determined by $\gamma$.  Let
$\Gamma_{\gamma}$ be the stabilizer in $\Gamma$ of $\gamma$, and we assume
that $\Gamma_{\gamma}$ is generated by $\gamma$; it is easily shown
that $\Gamma_{\gamma}$ is isomorphic to $\mathbb Z$.   Choose a realization of
$\Gamma$ in $\textrm{\rm PSL}_{2}(\mathbb R)$ such that $\gamma$ is diagonal.
Then the geodesic in $\mathbb H$ fixed by $\gamma$ is the line ${\cal L}_{0} =
\{\textrm{\rm Re}(z) = 0\}\cap \mathbb H$.  For any point $z \in M$, which we lift to
a point $z \in \mathbb H$, let $d_{\textrm{\rm hyp}}(z,{\cal L}_{0})$ denote the
geodesic distance from $z$ to ${\cal L}_{0}$.  With all this, we define the
\it hyperbolic counting function \rm as
$$
N_{\textrm{hyp};M,\gamma}(T;z) = \textrm{card }
\{ \eta \in \Gamma_{\gamma} \backslash \Gamma :
d_{\textrm{\rm hyp}}(\eta z, {\cal L}_{0}) < T\}\,.
$$
Equivalently, one can count the number of geodesic paths from $z \in M$
to the closed geodesic on $M$ in the homotopy class determined by $\gamma$; see
Figure 1.  By following the method of proof in Lemma 1.4 of \cite{JoLu}
(see also \cite{Lund}), which simply utilizes
elementary hyperbolic geometric considerations, we can establish the
following bound.  For any point $z \in M$ with injectivity radius $r$, and any
$u > T_{0} > r$, we have
\begin{equation}\label{countupperboundhyp}
N_{\textrm{hyp};M,\gamma}(u;z)\leq
 N_{\textrm{hyp};M,\gamma}(T_{0};z)+ \frac{\sinh^2
(\frac{u+r}{2})-\sinh^2 (\frac{T_0-r}{2})}{\sinh^2 (\frac{r}{2})}.
\end{equation}

\vskip .10in \noindent For the sake of completeness and
convenience of the reader, we now give a proof of
(\ref{countupperboundhyp}).

\vskip .10in Let $B_{z}(T)$ denote the hyperbolic ball of
hyperbolic radius $T$ centered at $z$.  Let $\{\eta_{k}\} \subset
\Gamma_{\gamma}\backslash \Gamma$ be a maximal collection of
elements such that $\eta_{k}z \in B_{z}(u) \setminus
B_{z}(T_{0})$.  Note that
$$
\bigcup_k B_{\eta_{k}z}(r) \subset B_{z}(u+r)\setminus
B_{z}(T_{0}-r),
$$
so then
$$
{\rm vol}_{\rm hyp}\left(\bigcup_k B_{\eta_{k}z}(r)\right) \leq
{\rm vol}_{\rm hyp}\left(B_{z}(u+r)\right) - {\rm vol}_{\rm
hyp}\left(B_{z}(T_{0}-r)\right),
$$
where ${\rm vol}_{\rm hyp}$ denotes the hyperbolic volume.  Since
$r$ is the injectivity radius at $z$, we then have
$$
\sum\limits_{k}{\rm vol}_{\rm hyp}\left( B_{\eta_{k}z}(r)\right)
\leq {\rm vol}_{\rm hyp}\left(B_{z}(u+r)\right) - {\rm vol}_{\rm
hyp}\left(B_{z}(T_{0}-r)\right).
$$
By computing the volume of geodesic balls, in $\mathbb{H}$, we
have that
$$
{\rm card}\{ \eta \in \Gamma_{\gamma}\backslash \Gamma : \eta z
\in B_{z}(u) \setminus B_{z}(T_0) \} \cdot  4\pi \sinh^2(r/2) \leq
4\pi \left( \sinh^2\left(\frac{u+r}{2} \right) -
\sinh^2\left(\frac{T_0 -r}{2}\right)\right).
$$
Since
$$
\begin{array}{rl}
{\rm card}\{ \eta \in \Gamma_{\gamma}\backslash \Gamma : \eta z
\in B_{z}(u) \setminus B_{z}(T_0) \} &= {\rm card}\{ \eta \in
\Gamma_{\gamma}\backslash \Gamma : \eta z
\in B_{z}(u)\}\\
\\ & - {\rm card}\{ \eta \in \Gamma_{\gamma}\backslash \Gamma : \eta z
\in B_{z}(T_0)\},
\end{array}
$$
we get the desired result, namely the bound
$$
{\rm card}\{ \eta \in \Gamma_{\gamma}\backslash \Gamma : \eta z
\in B_{z}(u) \} \leq {\rm card}\{ \eta \in
\Gamma_{\gamma}\backslash \Gamma : \eta z \in B_{z}(T_0)\} +
\frac{\sinh^2(\frac{u+r}{2}) - \sinh^2(\frac{T_0 -
r}{2})}{\sinh^2(\frac{r}{2})},
$$
thus completing the proof of (\ref{countupperboundhyp}).

\vskip .10in
Consider now a parabolic element $\gamma \in \Gamma$, which, by
conjugation in $\textrm{PSL}_{2}(\mathbb R)$, we may assume
$$
\gamma = \left(\begin{array}{cc} 1 & \omega \\ 0 & 1 \end{array} \right),
$$
where $\omega$ is referred to as the width of the cusp associated
to $\gamma$. Let $\Gamma_{\infty}$ denote the stabilizer in
$\Gamma$ of $\gamma$, and without loss of generality we may assume
that $\gamma$ generates $\Gamma_{\infty}$. Choose and fix any
point $z \in M$, which we lift to a point $z \in \mathbb H$.
Elementary considerations show that one can choose $y_{0} \in
\mathbb R$ sufficiently large so that $y_{0} > \textrm{\rm
Im}(\eta z)$ for all $\eta \in \Gamma$. Let ${\cal L}_{y_{0}}$ be
the horocyclic line in $\mathbb H$ defined by $\{\textrm{\rm
Im}(z) = y_{0}\}$.
 For any point $z \in M$, which we lift to
a point $z \in \mathbb H$, let $d_{\textrm{\rm hyp}}(z,{\cal L}_{y_{0}})$ denote the
geodesic distance from $z$ to ${\cal L}_{y_{0}}$.
With all this, we define the \it parabolic counting function \rm associated to $\gamma$ and $y_{0}$
to be
$$
N_{\textrm{par};M,P}(T;z,y_0) = \textrm{card }
\{ \eta \in \Gamma_{\infty} \backslash \Gamma : d_{\textrm{\rm hyp}}(\eta z, {\cal L}_{y_0}) < T\}.
$$

\vskip .10in
Observe that when defining the parabolic counting function, we needed to use the length
from $z$ to a
horocyclic line ${\cal L}_{y_{0}}$ since the cusp is at infinite distance.  Such considerations are
not necessary when defining the hyperbolic counting function.  Finally, as with (\ref{countupperboundhyp}),
the arguments from \cite{Lund} apply to yield the following bound.  For any point $z \in M$ with
injectivity radius $r$, and any $u > T_{0} > r$, we have
\begin{equation}\label{countupperboundpar}
N_{\textrm{par};M,P}(u;z,y_0)\leq
N_{\textrm{par};M,P}(T;z,y_0)+ \frac{\sinh^2
(\frac{u+r}{2})-\sinh^2 (\frac{T_0-r}{2})}{\sinh^2 (\frac{r}{2})}.
\end{equation}
The proof of (\ref{countupperboundpar}) is similar to the proof of
(\ref{countupperboundhyp}) given above.
\end{nn}

\begin{nn}\label{2.3}
\textbf{Parabolic Eisenstein series.}
By now, the study of parabolic Eisenstein series associated to a cusp
$P$ on a non-compact, finite volume hyperbolic Riemann surface $M$
is a classical aspect of mathematics (see, for example, \cite{Hejhal},
\cite{Iwaniec} or \cite{Kub}).  To recall, for any $z \in \mathbb H$
and $s \in \mathbb C$ with $\textrm{\rm Re}(s) > 1$, we define the
parabolic Eisenstein series $E_{\textrm{par};M,P}(s,z)$ by
\begin{equation}\label{paraeisenstein}
E_{\textrm{par};M,P}(s,z)= \omega^{-s}
\sum_{\eta \in \Gamma_{\infty} \backslash \Gamma} \left( \textrm{Im } \eta z\right)^s.
\end{equation}
It is standard in the mathematical literature to normalize cusps so that the width
$\omega$ is equal to one.  We will work slightly more generally and, as a result, include
the multiplicative factor of $\omega^{-s}$ in (\ref{paraeisenstein}).  For any point
$z \in \mathbb H$ and $y_{0} \in \mathbb R$ with $\textrm{\rm Im}(z) < y_{0}$, we have
that
$$
d_{\textrm{\rm hyp}}(z,{\cal L}_{y_{0}})
= \int_{\textrm{Im}(z)}^{y_0} \frac{dy}{y}
=\log\left(\frac{y_0}{\textrm{Im}(z)}\right)\,,
$$
so then
$$
\left( \textrm{Im}(z)\right)^s = y_{0}^{s}\exp
\left( -s \cdot d_{\textrm{\rm hyp}}(z,{\cal L}_{y_{0}}) \right)\,.
$$
With this observation, we can express the parabolic Eisenstein series (\ref{paraeisenstein})
as a Stieltjes integral, namely
\begin{equation}\label{paraint}
E_{\textrm{par};M,P}(s,z) = (y_0/\omega)^{s} \int_{0}^{\infty} e^{-su} dN_{\textrm{par};M,P}(u;z,y_0).
\end{equation}
Observe that the integral in (\ref{paraint}) depends on the choice of $y_{0}$ through
the parabolic counting function; however, after
multiplying by $y_{0}^{s}$ the product itself is independent of $y_{0}$.  As we will see,
one can use (\ref{countupperboundpar}) and (\ref{paraint}) to prove the well-known result that
(\ref{paraeisenstein}) converges uniformly and absolutely for $\textrm{\rm Re}(s) > 1$.
Though not needed in this article, we state, for the sake of completeness, the classical differential
equation satisfied by the parabolic Eisenstein series, which is the identity
$$
\Delta E_{\textrm{par};M,P} (s,z) = s(1-s) E_{\textrm{par};M,P}(s,z).
$$
\end{nn}

\begin{nn}\label{2.4}
\textbf{Hyperbolic Eisenstein series.}
Let $M = \Gamma \backslash \mathbb H$ be any finite volume, compact or non-compact,
hyperbolic Riemann surface, and let $\gamma$ be any hyperbolic element of $\Gamma$.
As in section \ref{2.2}, we assume that $\Gamma$ has been conjugated by an element
in $\textrm{\rm PSL}_{2}(\mathbb R)$ so that $\gamma$ is diagonal.  We will use
the change of coordinates $z = e^{\rho}e^{i\theta}$ and write $\theta(z) = \theta$.  With
this notation, we formally define the hyperbolic Eisenstein series
$E_{\textrm{hyp};M,\gamma}(s,z)$
by
\begin{equation}\label{hypeisenstein}
E_{\textrm{hyp};M,\gamma}(s,z) = \sum_{\eta \in \Gamma_{\gamma} \backslash \Gamma}
\left(\sin \theta(\eta z) \right)^s\,.
\end{equation}
The hyperbolic metric in the $(\rho, \theta)$ coordinates was given in section \ref{2.1},
from which one can easily show that
$$
d_{\textrm{\rm hyp}}(z, {\cal L}_{0}) = \vert\log(\csc \theta(z) +
\cot \theta(z) )\vert\,,
$$
which can be used to obtain the relation
$$
\sin (\theta (z)) \cdot \cosh (d_{\textrm{\rm hyp}}( z, {\cal L}_{0})) =1\,,
$$
so then we can write (\ref{hypeisenstein}) as
$$
E_{\textrm{hyp};M,\gamma}(s,z) = \sum_{\eta \in \Gamma_{\gamma} \backslash
\Gamma} \left(\cosh (d_{\textrm{\rm hyp}}(\eta z, {\cal L}_{0}))\right)^{-s}\,.
$$

We can express the hyperbolic Eisenstein series (\ref{hypeisenstein})
as a Stieltjes integral, namely
\begin{equation}\label{hypint}
E_{\textrm{hyp};M,\gamma}(s,z) = \int\limits_{0}^{\infty}(\cosh
u)^{-s}dN_{\textrm{hyp};M,\gamma}(u;z)\,.
\end{equation}
A by-product of the computations from section 4 is that by combining
(\ref{countupperboundhyp}) and (\ref{hypint}), we can show that the
series defining the hyperbolic Eisenstein series (\ref{hypeisenstein}) converges
uniformly and absolutely for $\textrm{\rm Re}(s) > 1$ (see also \cite{KM} and
\cite{Pippich}).  Also, using the computations from section \ref{2.1}, one
can easily verify the differential equation
\begin{equation}\label{hypeisendiff}
\Delta E_{\textrm{hyp};M,\gamma}(s,z) = s(1-s)E_{\textrm{hyp};M,\gamma}(s,z) +
s^{2}E_{\textrm{hyp};M,\gamma}(s+2,z)\,,
\end{equation}
which is given in \cite{KM}, \cite{Pippich} and \cite{Risager}
(Lemma 3.2).
\end{nn}

\begin{nn}\label{2.5}
\textbf{Degenerating families of Riemann surfaces.} The discussion
in this section is taken from \cite{JoLu2} and is repeated here
for the convenience of the reader.

\vskip .10in In \cite{JLu1} and \cite{JLu2} the authors gave a
construction of a degenerating family $M_{\ell}$ of either compact
or non-compact hyperbolic Riemann surfaces of finite volume. The
construction of the family $M_{\ell}$ allows one to define
unambiguously various notions such as the tracking of points
through degeneration and the idea of points not contained on the
degenerating geodesics. The reader is referred to these articles
for complete details, which will be assumed here.

\vskip .10in
Let $M_{\ell}$ be a degenerating family of connected,
hyperbolic Riemann surfaces with $p$ degenerating geodesics, with
$\ell$ denoting a $p$-tuple corresponding to the lengths of the degenerating geodesics.
To say that $\ell$ approaches zero means that the length of each degenerating geodesic is approaching zero.
 Although each $M_{\ell}$ is connected when $\ell > 0$, the limit
surface $M_{0}$ need not be connected and, indeed,
the number of cusps on $M_{0}$ is equal to the number of
cusps on $M_{\ell}$ plus $2p$.

\vskip .10in
For $\ell>0$, let ${\cal C}_{\ell}$ denote the hyperbolic infinite cylinder
with simple, closed geodesic of length $\ell$.
A convenient fundamental domain for ${\cal C}_{\ell}$ in $\mathbb H$ is
\begin{equation}\label{infinitecylinder}
\{r \exp(i\alpha): 1\leq r < \exp(\ell),\;\;0<\alpha<\pi\},
\end{equation}
with hyperbolic metric induced from $\mathbb H$ and uniformizing
group $\{\exp(k\ell):k \in \mathbb Z\}$ which acts on $\mathbb H$
by multiplication.  For any $\varepsilon > 0$, let ${\cal
C}_{\ell,\varepsilon}$ denote the symmetric submanifold of ${\cal
C}_{\ell}$ about the geodesic defined by $\gamma$ with total
volume equal to $\varepsilon$. A model for ${\cal
C}_{\ell,\varepsilon}$ in (\ref{infinitecylinder}) is obtained by
adding the restriction
$$
\cot^{-1}(\varepsilon/(2\ell)) <\alpha<\pi-\cot^{-1}(\varepsilon/(2\ell)).
$$
An easy calculation shows that the
length of each boundary component of ${\cal C}_{\ell,\varepsilon}$ is
$(\ell^2+\varepsilon^2 /4)^{1/2}$.  If $\varepsilon_{1} > \varepsilon_{0}$, then
the distance from the boundary of ${\cal C}_{\ell, \varepsilon_{1}}$ and
${\cal C}_{\ell, \varepsilon_{0}}$ can be shown to be
$$
d_{\textrm{\rm hyp}}(\partial {\cal C}_{\ell, \varepsilon_{0}}, \partial {\cal C}_{\ell, \varepsilon_{1}})=
\log \left((\varepsilon_{1} + \sqrt{(\varepsilon_{1}^{2}+4\ell^{2})})/
(\varepsilon_{0} + \sqrt{(\varepsilon_{0}^{2}+4\ell^{2})})\right).
$$

\vskip .10in From \cite{Randol} we have that for any $0 <
\varepsilon < 1/2$, the surface ${\cal C}_{\ell,\varepsilon}$
embeds isometrically into $M_{\ell}$.
 The surface $M_{0}$ contains $2p$ embedded copies of
${\cal C}_{0,\varepsilon}$ which is the limit of ${\cal
C}_{\ell,\varepsilon} \subset M_{\ell}$.  One can model ${\cal
C}_{0,\varepsilon}$ as two copies of a symmetric neighborhood of
the origin in the punctured unit disc with its complete hyperbolic
metric.   From \cite{Abikoff} we have that the family of
hyperbolic metrics converges uniformly on $M_{\ell} \setminus
{\cal C}_{\ell,\varepsilon}$.
\end{nn}

\begin{nn}\label{2.6}
{\bf A Stieltjes integral inequality.}
A key component in our analysis is an integral inequality for Stieltjes
integrals, which we quote from \cite{JK4} and, for the sake of completeness,
we state here.  Let $F$ be a real valued, smooth,
decreasing function defined for $u>0$ and let $g_{1},g_{2}$ be
real valued, non decreasing functions defined for $u\ge a>0$ and
satisfying $g_{1}(u)\le g_{2}(u)$ for $u\ge a$. Then, the
following inequality of Stieltjes integrals
$$
\int\limits_{a}^{\infty}F(u)\,dg_{1}(u)+F(a)\,g_{1}(a)\le
\int\limits_{a}^{\infty}F(u)\,dg_{2}(u)+F(a)\,g_{2}(a)
$$
holds, provided both integrals exist.
\end{nn}

\section{Convergence of counting functions}

In this section we will establish the limiting behavior of the
counting functions $N_{\textrm{par};M_{\ell},P}$ and
$N_{\textrm{hyp};M_{\ell},\gamma}$ on a degenerating family of
finite volume hyperbolic Riemann surfaces $M_{\ell}$.  For
simplicity, we will assume that $M_{\ell}$ has a single family of
degenerating geodesics; the more general situation is easily
obtained from the arguments presented here with only a slight
modification of notation.

\vskip .10in Throughout this article we make use of the following
fundamental result which we cite without proof from
\cite{Abikoff}, stated as Theorem 8, page 37.

\vskip .10in
\begin{nn}\label{3.1}
\textbf{Proposition.} {\it With notation as above, the hyperbolic
metrics on the degenerating family $M_{\ell}$ convergence to the
hyperbolic metric on $M_{0}$.  Furthermore, the convergence is
uniform on compact subsets of $M_{0}$ bounded away from the
developing cusps.}
\end{nn}

\vskip .10in We refer the reader to \cite{Abikoff} for a complete
proof of Proposition \ref{3.1}.  Building on this result, we
consider the convergence of the hyperbolic and parabolic counting
functions through degeneration.

\vskip .10in
\begin{nn}\label{3.2}
\textbf{Lemma.} {\it With notation as above, we have the following limits:}
\begin{list}{}
\item{a)} {\it If $\gamma$ does not correspond to a degenerating hyperbolic
element, then}
$$\lim\limits_{\ell \rightarrow 0}N_{\textrm{hyp};M_{\ell},\gamma}(T;z) =
N_{\textrm{hyp};M_{0},\gamma}(T;z);$$
\item{}
\item{b)} {\it For any cusp $P$, we have}
$$
\lim\limits_{\ell \rightarrow 0}N_{\textrm{par};M_{\ell},P}(T;z,y_0)
=N_{\textrm{par};M_{0},P}(T;z,y_0).
$$
\end{list}
{\it In all instances, the convergence is uniform on compact
subsets of $M_{0}$ bounded away from the developing cusps.}
\end{nn}

\proof Choose $\varepsilon_{1}$ sufficiently small so that the
point $z$ lies in $M_{\ell} \setminus {\cal
C}_{\ell,\varepsilon_{1}}$.  Now choose $\varepsilon_{0} <
\varepsilon_{1}$ so that the distance from the boundary of ${\cal
C}_{\ell,\varepsilon_{0}}$ to the boundary of ${\cal
C}_{\ell,\varepsilon_{1}}$ is greater than $T$. Clearly, any
geodesic path from $z$ to the non-pinching geodesic $\gamma$ with
length bounded by $T$ necessarily lies entirely in the $M_{\ell}
\setminus {\cal C}_{\ell,\varepsilon_{0}}$. From Proposition
\ref{3.1}, we know that the family of hyperbolic metrics converge
uniformly away on $M_{\ell} \setminus {\cal
C}_{\ell,\varepsilon_{0}}$, which proves part (a).

\vskip .10in The convergence statement asserted in (b) follows from a similar argument.
\hfill $\Box$

\begin{nn}\label{3.3}
\textbf{Lemma.} {\it With notation as above, let}
$$
g(y_0,\ell) = \log \left(\frac{\omega}{y_0 \ell} +
\sqrt{\left(\frac{\omega}{y_0 \ell}\right)^2+1}\right)\,
$$
{\it and let $\gamma$ correspond to a degenerating hyperbolic element.}
\begin{list}{}
\item{a)} {\it Assume $\varepsilon > 0$ is sufficiently small so that
${\cal C}_{\ell,\varepsilon}$ is embedded in $M_{\ell}$, and, for $z \in M_{\ell}$, let
${\cal C}_{\ell,\varepsilon}^{z}$ denote the half of ${\cal C}_{\ell,\varepsilon}$ closest to $z$.
Define}
\begin{equation}\label{boundarycount}
N_{\textrm{hyp};M_{\ell},\partial {\cal C}^{z}_{\ell,\varepsilon}}(T;z) = \textrm{card }
\{ \eta \in \Gamma_{\gamma} \backslash \Gamma :
d_{\textrm{\rm hyp}}(\eta z, \partial {\cal C}^{z}_{\ell,\varepsilon}) < T\}\,.
\end{equation}
{\it Then, for any $T > 0$, we have}
$$
 N_{\textrm{hyp};M_{\ell},\gamma}(T + g(y_{0},\ell);z)
= N_{\textrm{hyp};M_{\ell},\partial {\cal C}^{z}_{\ell,\varepsilon}}(T;z)\,.
$$
\item{b)} {\it For any fixed $T > 0$, we have that}
$$
\lim\limits_{\ell \rightarrow 0} N_{\textrm{hyp};M_{\ell},\gamma}(T+g(y_{0},\ell);z)
= N_{\textrm{par};M_{0},P}(T;z,y_0).
$$
\end{list}
\end{nn}

\begin{figure}
\begin{center}
\includegraphics[width=0.6\textwidth]{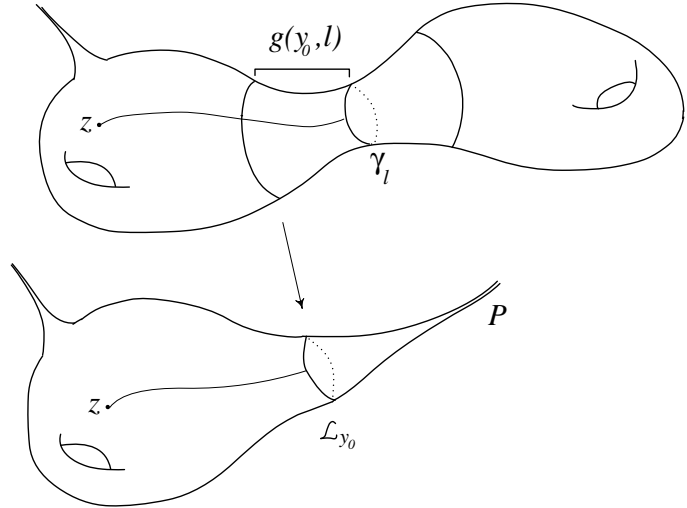}
\caption{The Riemann surfaces $M_{\ell}$ and $M_0$}
\end{center}
\end{figure}

\proof
For fixed $\ell$,
let us identify $M_{\ell}$ with a (Ford) fundamental domain in $\mathbb H$ such
that the lift of the pinching geodesic $\gamma$ lies along the line
$\textrm{\rm Re}(z) = 0$.  Then, the boundary $\partial {\cal C}^{z}_{\ell,\varepsilon}$
of ${\cal C}_{\ell,\varepsilon}$ lies along a ray $\theta(z) = \textrm{\rm constant}$.
The curve $\partial{\cal C}_{\ell,\varepsilon}^z$ is orthogonal to the geodesics which transverse the
sub-cylinder ${\cal C}_{\ell, \varepsilon}^{z}$, so $\partial{\cal C}_{\ell,\varepsilon}^{z}$
converges to a path on $M_{0}$ which is perpendicular to the geodesics which transverse a
neighborhood of the cusp, meaning $\partial{\cal C}_{\ell, \varepsilon}$ converges to a horocyclic path
${\cal L}_{y_{0}}$ in a neighborhood of the new cusp; see Figure 2.
The area of ${\cal C}_{\ell, \varepsilon}^{z}$
equals the area on $M_{0}$ above ${\cal L}_{y_{0}}$ by the choice of $\varepsilon = \varepsilon (\ell)$.
By direct computation, we have
$$
\textrm{Area}(\textrm{\rm region on $M_{0}$ above} \,\,{\cal{L}}_{y_{0}})
 = \int_{y_0}^{\infty} \int_{0}^{\omega} \frac{dxdy}{y^2} =
\frac{\omega}{y_0},
$$
thus we get the relation $\varepsilon = 2\omega/y_{0}$.
We define $g(y_{0},\ell)$ to be the distance from $\partial {\cal C}_{\ell,\varepsilon}^{z}$
to the geodesic in ${\cal C}_{\ell,\varepsilon}^{z}$.  Using the computations from section \ref{2.1}
and section \ref{2.5}, we then have that
$$
g(y_{0},\ell) = \int\limits_{\cot^{-1}(\varepsilon/2\ell)}^{\pi/2} \frac{d\theta}{\sin \theta}\,,
$$
which is easily evaluated, arriving at the claimed result.

Choose any $\eta \in \Gamma_{\gamma}\backslash \Gamma$.  By
choosing the appropriate coset representative, we may assume that
$\eta z$ lies in the fundamental domain for $\Gamma_{\ell}$ from
section \ref{2.5}, meaning $0 \leq \log \vert \eta z \vert <
\ell$.  It is immediate that the geodesic path from $\eta z$ to
the $\{\textrm{\rm Re}(z) = 0\} \cap \mathbb H$ lies along the
path $\rho = \textrm{\rm constant}$, which then is seen to be
orthogonal to each ray $\theta = \textrm{\rm constant}$.
Therefore, we have that
\begin{equation}\label{triangleequality}
d_{\textrm{\rm hyp}}(\eta z, \gamma) =
d_{\textrm{\rm hyp}}(\eta z, \partial {\cal C}_{\ell,\varepsilon}^{z}) +
d_{\textrm{\rm hyp}}(\partial {\cal C}_{\ell,\varepsilon}^{z}, \gamma)
= d_{\textrm{\rm hyp}}(\eta z, \partial {\cal C}_{\ell,\varepsilon}^{z}) + g(y_{0},\ell).
\end{equation}
From (\ref{triangleequality}) it is clear that
$d_{\textrm{\rm hyp}}(\eta z, \gamma) < T + g(y_{0},\ell)$ if and only if
$d_{\textrm{\rm hyp}}(\eta z, \partial {\cal C}_{\ell,\varepsilon}^{z}) < T$, which
completes the proof of part (a).

Part (b) follows from combining part (a) with the convergence of
the hyperbolic metric on $M_{\ell}$ away from the developing cusps
to the hyperbolic metric on $M_{0}$, as stated in Proposition
\ref{3.1}. \hfill $\Box$

\begin{nn}\label{3.4}
\textbf{Remark.} As discussed after the proof of the main theorem,
there are two cases one needs to consider in part (b) of Lemma
\ref{3.3}:  When the degenerating geodesic is separating, and when
the degenerating geodesic is non-separating.  If $\gamma$ is
separating, then the statement of (b) holds without any liberty in
the notation.  If $\gamma$ is non-separating, however, one needs
to take into account that geodesic lengths from $z$ to $\gamma$
enter the cylinder about the pinching geodesic from the two
different sides.  The proof of (b) immediately extends to show
that in the non-separating case the right-hand-side is actually
the sum of two parabolic counting functions corresponding to the
two newly formed cusps.  With this noted, we choose to use the
statement in (b) with its slight abuse in both the separating and
non-separating cases in order to prevent burdensome notation.
\end{nn}

\begin{nn}\label{3.5}
\textbf{Remark.} As one can see, the convergence of the counting
functions in Lemma \ref{3.2} follows directly from the convergence
of the hyperbolic metrics away from the developing cusps, as
stated in Proposition \ref{3.1}.  In Lemma \ref{3.3}, we have the
added feature that the hyperbolic counting function involves the
distances from the orbits of a point $z$ to the geodesic
corresponding to the hyperbolic element, but the parabolic
counting function involves distances to a chosen horocycle.  The
distances to the geodesic associated to a degenerating hyperbolic
element are growing without bound; however, Lemma \ref{3.3} can be
viewed as establishing a type of ``regularized convergence''.  To
be more specific, observe that the function $g(y_{0},\ell)$
depends solely on $y_{0}$ and $\ell$, and no other aspect of the
family $M_{\ell}$.  With this in mind, Lemma \ref{3.3} states that
if we ``regularized'' the counting functions $N_{{\rm hyp};
M_{\ell}, \gamma}$ by introducing the factor $g(y_{0},\ell)$, one
then has convergence of the counting functions.

\end{nn}

\section{Convergence of Eisenstein series}

In this section, we prove the Main Theorem.  In brief, our proof
uses the convergence of the counting functions for fixed $T$
(Lemma \ref{3.2} and Lemma \ref{3.3}), the uniform bounds for the
counting functions (section \ref{2.2}) and the Stieltjes integral
inequality (section \ref{2.6}).  As in section 3, we present the
arguments in the setting of a single degenerating hyperbolic
element $\gamma$ whose geodesic has length $\ell$; in order to
consider the general situation where there are a number of
degenerating geodesics, one simply needs notational changes.

\begin{nn}\label{4.1}
\textbf{Proof of the main theorem.} {\it Proof of part (i):} For any $T_{0} > 0$, write
\begin{equation}\label{twointegrals}
E_{\textrm{hyp};M_{\ell},\gamma}(s,z) =
\int_{0}^{T_0} (\cosh u)^{-s} d N_{\textrm{hyp};M_{\ell},\gamma}(u;z)
+ \int_{T_0}^{\infty} (\cosh u)^{-s} d N_{\textrm{hyp};M_{\ell},\gamma}(u;z) .
\end{equation}
Choose any $\delta > 0$ and restrict $s \in \mathbb C$ to the half-plane
$\textrm{\rm Re}(s) \geq 1 + \delta$ for some fixed $\delta > 0$.  Trivially, we have
that
$$
\left|\int_{T_0}^{\infty} (\cosh u)^{-s} dN_{\textrm{hyp};M_{\ell},\gamma}(u;z) \right|
\leq
\int_{T_0}^{\infty} (\cosh u)^{-(1+\delta)} dN_{\textrm{hyp};M_{\ell},\gamma}(u;z).
$$
We now establish the following bound:  Given any $\varepsilon > 0$,
there is a $T_{0} = T_{0}(\varepsilon, \delta, r)$, where $r$ is the injectivity radius
at $z$, such that for each $\ell \geq 0$, we have
\begin{equation}\label{tailbound}
\int_{T_0}^{\infty} (\cosh u)^{-(1+\delta)} dN_{\textrm{hyp};M_{\ell},\gamma}(u;z) <
\epsilon\,.
\end{equation}
The verification of (\ref{tailbound}) follows the proof of Lemma 1.4 from \cite{JoLu},
which we repeat here.  In the notation of section 2.6,
let $F(u)=(\cosh u)^{-(1+\delta)}$, which evidentally is real-valued, smooth, and decreasing.
For $u > T_{0}$, we let
\begin{eqnarray*}
g_1(u) &=&  N_{\textrm{hyp};M_{\ell},\gamma}(u;z)\\
g_2(u) &=& N_{\textrm{hyp};M_{\ell},\gamma}(T_{0};z)+ \frac{\sinh^2
(\frac{u+r}{2})-\sinh^2 (\frac{T_0-r}{2})}{\sinh^2 (\frac{r}{2})}.
\end{eqnarray*}
As stated in section 2.2, we have that $g_1(u) \le g_2(u)$,
and both $g_{1}$ and $g_{2}$ are real-valued and non-decreasing for
$u \ge T_0>0$.  With all this, the Stieltjes integral inequality
from section 2.6 yields the bound
\begin{eqnarray*}
\int_{T_0}^{\infty} (\cosh u)^{-(1+\delta)}
dN_{\textrm{hyp};M_{\ell},\gamma}(u;z) & \le& \int_{T_0}^{\infty}
(\cosh u)^{-(1+\delta)} dg_2(u) \\ &  + & (\cosh
T_0)^{-(1+\delta)} \Bigg\{ \frac{\sinh^2 (\frac{T_0+r}{2})-\sinh^2
(\frac{T_0-r}{2})}{\sinh^2 (\frac{r}{2})} \Bigg\}.
\end{eqnarray*}
Elementary calculations and trigonometric identities imply that
$$
dg_2(u) = \frac{\sinh (u+r)}{2\sinh^2 (\frac{r}{2})} du
$$
and
$$
\sinh^2 (\frac{T_0+r}{2})-\sinh^2 (\frac{T_0-r}{2}) =\sinh r \sinh
T_0\,.
$$
Using the trivial bounds $\sinh u \leq e^{u}/2$ and $\cosh u \geq e^{u}/2$, we
then obtain the estimates
\begin{equation}\label{upperbound}
\begin{array}{lll}\displaystyle
\int_{T_0}^{\infty} (\cosh u)^{-(1+\delta)} dg_1(u) &\displaystyle \le& \displaystyle \frac{1}{2
\sinh^2(\frac{r}{2})}\int_{T_0}^{\infty} (\cosh u)^{-(1+\delta)} \sinh(u+r)
du
\\[5mm]&\displaystyle  + &\displaystyle
  (\cosh T_0)^{-(1+\delta)} \frac{\sinh r \sinh T_0}{\sinh^2 (\frac{r}{2})} \\[5mm]
&\displaystyle \le &\displaystyle  \frac{2^{\delta}e^r}{\sinh^2(\frac{r}{2})}\int_{T_0}^{\infty}
e^{-\delta\cdot u}du
+
(\cosh T_0)^{-(1+\delta)} \frac{\sinh r \sinh T_{0}}{\sinh^2 (\frac{r}{2})}\\[5mm]
&\displaystyle \le&\displaystyle
 e^{-\delta\cdot T_{0}} \left(\frac{2^{\delta}e^r}{\delta \sinh^2(\frac{r}{2})}
+ \frac{2^{\delta}\sinh r}{\sinh^2 (\frac{r}{2})} \right) \,,
\end{array}
\end{equation}
which clearly can be made smaller than any $\varepsilon > 0$, namely, by taking
\begin{equation}\label{lowerbound}
T_{0} \geq \frac{1}{\delta}\left(-\log \varepsilon + \log
\left(\frac{2^{\delta}e^r}{\delta \sinh^2(\frac{r}{2})}
+ \frac{2^{\delta}\sinh r}{\sinh^2 (\frac{r}{2})}\right)\right)\,.
\end{equation}
Therefore, we have proved the bound asserted in (\ref{tailbound}).

\vskip .10in In addition to (\ref{lowerbound}) let us assume, for
convenience, that $T_{0}$ is a point of continuity of
$N_{\textrm{hyp};M_{0},\gamma}(T;z)$, meaning there is no geodesic
path from $z$ to $\gamma$ on $M_{0}$ with length equal to $T_{0}$.
Then, with $T_{0}$ chosen, there is an integer $N$ and an
$\ell_{0}$ sufficiently small such that for $\ell < \ell_{0}$, we
have $N=N_{\textrm{hyp};M_{\ell},\gamma}(T_{0};z)  =
N_{\textrm{hyp};M_{0},\gamma}(T_{0};z)$.  Let $\{d_{k,M_{\ell}}\}
\subset [0,T_{0}]$ be the set of lengths on $M_{\ell}$ such that
for any $\eta > 0$ we have
$$
N_{\textrm{hyp};M_{\ell},\gamma}(d_{k,M_{\ell}}-\eta;z) <
N_{\textrm{hyp};M_{\ell},\gamma}(d_{k,M_{\ell}}+\eta;z).
$$
For simplicity, we count the elements in the set
$\{d_{k,M_{\ell}}\}$ with multiplicities so that we have
$$
\int_{0}^{T_0} (\cosh u)^{-s} d
N_{\textrm{hyp};M_{\ell},\gamma}(u;z) =
\sum\limits_{k=1}^{N}\left(\cosh d_{k,M_{\ell}}\right)^{-s}.
$$
With this, we can write
\begin{eqnarray*}\label{firstterms}
\int_{0}^{T_0} (\cosh u)^{-s} d N_{\textrm{hyp};M_{\ell},\gamma}(u;z) &-&
\int_{0}^{T_0} (\cosh u)^{-s} d N_{\textrm{hyp};M_{0},\gamma}(u;z) \\
&=& \sum\limits_{k=1}^{N} \left[(\cosh d_{k,M_{\ell}})^{-s} -  (\cosh d_{k,M_{0}})^{-s}
\right]\, .
\end{eqnarray*}
Observe now that the function $(\cosh u)^{-s}$ is uniformly
continuous and absolutely continuous on $[0,T_{0}]$. By Lemma
\ref{3.2}, which we apply for all $T < T_{0}$, there is an
$\ell_0'$ such that for $\ell<\ell_0'$ we have
$$
|d_{k,M_{\ell}}-d_{k, M_0}| < \frac{\delta}{N} \, \textrm{ for all
$k$,}
$$
so then, $\sum_{k=1}^{N}|d_{k,M_{\ell}}-d_{k, M_0}| < \delta$. By
the absolute continuity of $(\cosh u)^{-s}$ on $[0,T_0]$ we arrive
at the bound
\begin{equation}\label{firsttermbound}
\left| \sum\limits_{k=1}^{N} \left[(\cosh d_{k,M_{\ell}})^{-s} -  (\cosh d_{k,M_{0}})^{-s}
\right] \right| \le \sum\limits_{k=1}^{N} \left|(\cosh d_{k,M_{\ell}})^{-s} -  (\cosh d_{k,M_{0}})^{-s}
\right| < \varepsilon \,.
\end{equation}

To put all this together, let us write
\begin{equation}\label{triangle}
\begin{array}{l}
\displaystyle \left|
E_{\textrm{hyp};M_{\ell},\gamma}(s,z) - E_{\textrm{hyp};M_{0},\gamma}(s,z) \right| \le \\[5mm]
\displaystyle \hskip 1.5in \left| \int_{0}^{T_0} (\cosh u)^{-s} d N_{\textrm{hyp};M_{\ell},\gamma}(u;z) -
\int_{0}^{T_0} (\cosh u)^{-s} d N_{\textrm{hyp};M_{0},\gamma}(u;z) \right| \\[5mm]
\hskip 1.5in \displaystyle  + \,\,\,\,
 \left| \int_{T_0}^{\infty} (\cosh u)^{-s} d N_{\textrm{hyp};M_{\ell},\gamma}(u;z) \right| \\[5mm]
\hskip 1.5in \displaystyle + \,\,\,\,
\left| \int_{T_0}^{\infty} (\cosh u)^{-s} d N_{\textrm{hyp};M_{0},\gamma}(u;z) \right|\,.
\end{array}
\end{equation}

The second and third terms on the right hand side are arbitrarily small by taking $T_{0}$ as in
(\ref{lowerbound}), and the first term on the right hand side is arbitrarily small
by (\ref{firsttermbound}). With all this, the proof of part (i) of the main theorem is complete.

\begin{nn}\label{4.2}
\textbf{Remark.} The referee has proposed the following alternate
proof of (\ref{tailbound}).  For any given $\ell$, on the geodesic
$\gamma$ there are finitely many points $w_{\ell,j} , j = 1,
\ldots ,K$ which partition $\gamma$ into segments of length $<
\delta_1$.  Since $\gamma$ is not a pinching geodesic, we can take
$\gamma$, as well as the partitioning points, as lying in a subset
of $M_{\ell}$ which is bounded away from the developing cusps. For
any $\eta \in \Gamma_{\gamma} \backslash \Gamma$ with $d(\eta
z,\mathcal{L}_0) < T$ on $M_{\ell}$, let $w_{\ell}$ be the point
on $\gamma$ such that $d(\eta z,w_{\ell}) = d(\eta
z,\mathcal{L}_0)$.  Using that $w_{\ell}$ is within distance
$\delta_1$ from some $w_{\ell,j}$, the triangle inequality gives
the bound
$$
d(\eta z, w_{\ell,j}) \le d(\eta z, w_{\ell}) + \delta_1 \le T +
\delta_1.
$$
If we let $N_{\Gamma}(z,w,t)$ denote the counting function for the
groups elements that move $z$ within distance $t$ from $w$, we
then arrive at the inequality (hyperbolic lattice counting).
$$
N_{\textrm{hyp};M_{\ell},\gamma}(T;z) \le \sum_{j}
N_{\Gamma}(z,w_{\ell,j},T+\delta_1).
$$
Using hyperbolic volume considerations, one trivially shows that
$N_{\Gamma}(z,w_{\ell,j},T+\delta_1)$ is bounded by
$O(e^{T+\delta_1})$, and the bound is uniform for the $w_{\ell,j}$
contained in a compact set.  Therefore, we can write
$$
N_{\textrm{hyp};M_{\ell},\gamma}(T;z) \ll e^{T+\delta_1}.
$$
Returning to (\ref{tailbound}), one can integrate by parts to get
\begin{eqnarray*}
\int_{T_0}^{\infty} (\cosh u)^{-1-\delta} d
N_{\textrm{hyp};M_{\ell},\gamma}(u;z)
&=& \left[ (\cosh u)^{-1-\delta} N_{\textrm{hyp};M_{\ell},\gamma}(u;z) \right]_{T_0}^{\infty}\\
&+&(1+\delta)\int_{T_0}^{\infty} (\cosh u)^{-2-\delta} \sinh u
N_{\textrm{hyp};M_{\ell},\gamma}(u;z)du.
\end{eqnarray*}
The discussion above implies that
$N_{\textrm{hyp};M_{\ell},\gamma}(u;z)$ is bounded by
$O\left(e^{u+\delta_1}\right)$ independently of $\ell$. If we take
$\delta_1 < \delta$, then we can easily choose $T_0$ with the
required property that the original integral is $< \epsilon$.
Indeed, the first term is $O(e^{(-\delta+\delta_1)T_0})$ and the
same applies to the second since $(\cosh u)^{-2-\delta} \sinh u =
O(e^{(-1-\delta)u})$.

\vskip .10in As noted by the referee, an important aspect of the
above argument is that one only needs the rough order of growth of
$N_{\textrm{hyp};M_{\ell},\gamma}(u;z)$, i.e. the injectivity
radius plays no role in the formula.
\end{nn}

\vskip .10in
{\it Proof of part (ii):}
The proof of part (ii) follows the pattern set in the proof of part (i).  The only difference
is that one is considering the function $F(u) = e^{-su}$, rather than $F(u) = (\cosh u)^{-s}$.
For the integral over $[T_{0},\infty)$, the essential feature from $F$ to be used is that
$\vert F(u)\vert e^{u}$ is integrable.  For the integral over $[0,T_{0}]$, one needs $F$ to be
uniformly and absolutely continuous.

\vskip .10in
{\it Proof of part (iii):}
We proceed as in the proof of parts (i) and (ii) with a few
slight modifications.  To begin, we write
\begin{equation}\label{twointegrals2}
E_{\textrm{hyp};M_{\ell},\gamma}(s,z) =
\int_{0}^{T_0+g(y_{0},\ell)} (\cosh u)^{-s} d N_{\textrm{hyp};M_{\ell},\gamma}(u;z)
+ \int_{T_0+g(y_{0},\ell)}^{\infty} (\cosh u)^{-s} d N_{\textrm{hyp};M_{\ell},\gamma}(u;z)\,,
\end{equation}
where $g(y_{0},\ell)$ is given in Lemma \ref{3.3}.  We shall
multiply both sides of (\ref{twointegrals2}) by
$2^{-s}e^{sg(y_{0},\ell)}$ and let $\ell$ approach zero. For the
integral over $[T_{0}+g(y_{0},\ell),\infty)$, we first use part
(a) of Lemma \ref{3.3} to write
$$
\int_{T_0+g(y_{0},\ell)}^{\infty} (\cosh u)^{-s} d
N_{\textrm{hyp};M_{\ell},\gamma}(u;z) = \int_{T_0}^{\infty} (\cosh
(u+g(y_{0},\ell)))^{-s} d N_{\textrm{hyp};M_{\ell},\partial {\cal
C}_{\ell,\varepsilon}}(u;z)\,.
$$
The geometric argument from \cite{JoLu} and \cite{Lund}
which produced (\ref{countupperboundhyp}) and (\ref{countupperboundpar})
immediately extends to give the bound
$$
N_{\textrm{hyp};M_{\ell},\partial {\cal C}_{\ell,\varepsilon}}(u;z)\leq
N_{\textrm{hyp};M_{\ell},\partial {\cal C}_{\ell,\varepsilon}}(T_{0};z)+ \frac{\sinh^2
(\frac{u+r}{2})-\sinh^2 (\frac{T_{0}-r}{2})}{\sinh^2 (\frac{r}{2})}\,,
$$
for $u > T_{0} > r$ where, as before, $r$ is the injectivity radius of $M_{\ell}$ at $z$.
Following the computations in (\ref{upperbound}), we arrive at the estimate
\begin{equation}\label{hyppartail}
\left|2^{-s}e^{sg(y_0,\ell)} \int_{T_0+g(y_0,\ell)}^{\infty}
(\cosh u)^{-s} d N_{\textrm{hyp};M_{\ell},\gamma}(u;z) \right| \le
e^{-\delta\cdot T_{0}}
 \left(\frac{e^{r}}{\sinh^{2}(r/2)} \right)\,,
\end{equation}
where we have written $\textrm{\rm Re}(s) = 1 +\delta$. By
choosing
$$
T_{0} \geq \frac{1}{\delta}\left(-\log \varepsilon + \log
\left(\frac{e^{r}}{\sinh^{2}(r/2)}\right)\right),
$$
we have that the upper bound in
(\ref{hyppartail}) is less than $\varepsilon$.

For the first integral in (\ref{twointegrals2}), we begin by writing
$$
\int_{0}^{T_0+g(y_{0},\ell)} (\cosh u)^{-s} d
N_{\textrm{hyp};M_{\ell},\gamma}(u;z) = \int_{0}^{T_0}(\cosh
(u+g(y_{0},\ell)))^{-s} d N_{\textrm{hyp};M_{\ell},\partial {\cal
C}_{\ell,\varepsilon}}(u;z)\,.
$$
Also, we observe the following elementary result:
For fixed $x > 0$ and $s \in \mathbb C$ with $\textrm{\rm Re}(s) > 0$, we have
\begin{equation}\label{hypparlimit}
\lim\limits_{r\rightarrow \infty} 2^{-s}e^{rs}(\cosh (x+r))^{-s} = e^{-sx}\,.
\end{equation}
Furthermore, the limit (\ref{hypparlimit}) is uniform for all
$x>0$ and $\textrm{\rm Re}(s) \geq 1 + \delta$. Let
$f(s,\ell)=2^{-s}y_0^s e^{sg(y_0,\ell)}$. Then, by Lemma \ref{3.3}
and the argument yielding (\ref{firsttermbound}), we have, for any
$T_{0}$ as in (\ref{lowerbound}), the limit
\begin{equation}\label{hypparfirst}
\lim\limits_{\ell \rightarrow 0}
f(s,\ell) \int_{0}^{T_0+g(y_0,\ell) } (\cosh (u))^{-s} dN_{\textrm{hyp};M_{\ell},\gamma}(u;z)=
 y_0^s \int_{0}^{T_{0}} e^{-su} dN_{\textrm{par};M_{0},P}(u;z,y_0)\,.
\end{equation}

We now use (\ref{hyppartail}) and (\ref{hypparfirst}) and the triangle inequality, as
in (\ref{triangle}), in order to prove
\begin{equation}\label{partiii}
\lim\limits_{\ell \rightarrow 0} f(s,\ell) E_{{\rm
hyp};M_{\ell},\gamma}(s,z) = \omega^{s} E_{{\rm par};M_{0},P}(s,z)\,.
\end{equation}
To complete the proof of part (iii), it remains to evaluate $f(s,\ell)$.

\vskip .10in {\it Evaluation of $f(s,\ell)$.} As shown in the
proof of Lemma \ref{3.3}, we have
$$
g(y_0,\ell) = \log \Bigg( \frac{\omega}{y_0 \ell} + \sqrt{\Big(\frac{\omega}{y_0
\ell}\Big)^2+1}\Bigg),
$$
from which we immediately derive the relation
\begin{equation}\label{asymptotics}
f(s,\ell)=2^{-s}y_0^s \Bigg( \frac{\omega}{y_0 \ell} +
\sqrt{\Big(\frac{\omega}{y_0 \ell}\Big)^2+1}\Bigg)^s = (\omega/\ell)^{s} +
o\left( (\omega/\ell)^{s}\right)\,\,\,\,\,\textrm{\rm as} \,\,\,\,\,
\ell \rightarrow 0.
\end{equation}
Substituting (\ref{asymptotics}) into (\ref{partiii}), then multiplying both
sides by $\omega^{-s}$, completes the proof of part (iii) of our Main Theorem.
\end{nn}

\begin{nn}\label{4.3}
\textbf{Remark.} In the setting of part (iii) of our Main Theorem,
consider the differential equation satisfied by $E_{{\rm
hyp};M_{\ell},\gamma}(s,z)$ which, after multiplying by
$\ell^{-s}$, is the identity
\begin{equation}\label{diffeq}
\Delta (\ell^{-s} E_{\textrm{hyp};M,\gamma}(s,z)) =
s(1-s)(\ell^{-s} E_{\textrm{hyp};M,\gamma}(s,z)) +
(s\ell)^{2}(\ell^{-s-2}E_{\textrm{hyp};M,\gamma}(s+2,z))\,.
\end{equation}
By part (iii) of our Main Theorem, we have that
$$
\lim\limits_{\ell \rightarrow 0}\left(s(1-s)(\ell^{-s} E_{\textrm{hyp};M,\gamma}(s,z)) +
s^{2}(\ell^{-s-2}E_{\textrm{hyp};M,\gamma}(s+2,z))\cdot \ell^{2}\right) = s(1-s) E_{{\rm par};M_{0},P}(s,z)\,,
$$
for all $\textrm{\rm Re}(s) > 1$ and $z$ bounded away from the developing cusps.  The point here
is that the second term on the right hand side of (\ref{diffeq}) vanishes through
degeneration.  Heuristically, this shows that in the setting of part (iii),
the differential equation for the hyperbolic Eisenstein series
limits to the differential equation for the parabolic Eisenstein series.
\end{nn}

\begin{nn}\label{4.4}
\textbf{Remark.}
In the definition of the parabolic Eisenstein series (\ref{paraeisenstein}) we included
a multiplicative factor of $\omega^{-s}$.  Analogously, we could have included a factor
of $\ell_{\gamma}^{-s}$ in the definition of the hyperbolic Eisenstein series (\ref{hypeisenstein}).
Let us use the term \it adjusted hyperbolic Eisenstein series \rm to denote the hyperbolic
Eisenstein series from (\ref{hypeisenstein}) multiplied by $\ell_{\gamma}^{-s}$.
With this factor, then part (iii) of the Main Theorem states that the adjusted hyperbolic
Eisenstein series associated to the degenerating hyperbolic element converges to the parabolic
Eisenstein series of the newly formed cusp.  In addition, the adjusted hyperbolic Eisenstein
series will satisfy an equation similar to (\ref{hypeisendiff}), where the second
term has the multiplicative factor of $(s\ell_{\gamma})^{2}$, as in (\ref{diffeq}).
\end{nn}

\begin{nn}\label{4.5}
\textbf{Remark.} The concept of an Eisenstein series associated to
an elliptic element of $\Gamma$ was first defined in \cite{JK11}
and has been studied in \cite{Pippich}.  At this time, A. Pippich
is continuing her systematic investigation of elliptic Eisenstein
series, which, almost certainly, will include convergence results
as in the present paper when considering a sequence of
elliptically degenerating Riemann surfaces.
\end{nn}

\begin{nn}\label{4.6}
\textbf{Remark.} After completion and initial review of this
article, Gautam Chinta called our attention to the article
\cite{Falliero} where the author establishes the Main Theorem
using different techniques.  The advantage of our approach is the
introduction of counting function techniques when studying
Eisenstein series, both parabolic and hyperbolic, thus reducing
the main theorem to convergence questions associated to the
various counting functions.
\end{nn}

\vskip .10in \textbf{Acknowledgements.}  The authors thank Anupam
Bhatnagar, James Cogdell, Jozek Dodziuk, and James Harlacher for
numerous helpful conversations.  In particular, we thank Cogdell
for several comments concerning exposition as well as the contents
of Remark \ref{4.4}.  We also thank the referee for providing an
extensive list of comments which enhanced the presentation and
clarity of the paper. The second named author (J.J.) acknowledges
support from NSF grants and several PSC-CUNY awards.

\vskip .40in  \vspace{5mm}\noindent
Daniel Garbin \\
Mathematics Ph.D. Program \\
The Graduate Center of CUNY \\
365 Fifth Avenue\\
New York, NY
U.S.A. \\
e-mail: daniel$\_$garbin@yahoo.com

\vspace{5mm} \noindent
Jay Jorgenson \\
Department of Mathematics \\
The City College of New York \\
Convent Avenue at 138th Street \\
New York, NY 10031
U.S.A. \\
e-mail: jjorgenson@mindspring.com

\vspace{5mm}\noindent
Michael Munn \\
Mathematics Ph.D. Program \\
The Graduate Center of CUNY \\
365 Fifth Avenue\\
New York, NY
U.S.A. \\
e-mail: mikemunn@gmail.com


\begin{thebibliography}{99}


\bibitem{Abikoff}
Abikoff, W.: Degenerating families of Riemann surfaces. \it Annals
of Math. \bf 105 \rm (1977) 29--44.

%\bibitem{Chavel}
%Chavel, Isaac
%Eigenvalues in Riemannian geometry.
%Pure and Applied Mathematics, 115.
%Academic Press, Inc., Orlando, FL, 1984. xiv+362 pp.

\bibitem{Fay}
Fay, J.: Theta functions on Riemann surfaces. Lecture Notes in
Mathematics \bf 352\rm. Springer-Verlag, Berlin-New York, (1973)
iv+137 pp.

\bibitem{Falliero}
Falliero, T.: D\'eg\'en\'erescence de s\'eries d'Eisenstein
hyperboliques. \it Math. Ann. \bf 339 \rm (2007) 341--375.

\bibitem{Hejhal}
Hejhal, D.: The Selberg trace formula for ${\rm PSL}(2,\,R)$. Vol.
2. Lecture Notes in Mathematics \bf 1001\rm. Springer-Verlag,
Berlin, (1983) viii+806 pp.

\bibitem{HejhalAMS}
Hejhal, D.:  Regular $b$-groups, degenerating Riemann surfaces,
and spectral theory, Memoirs of the AMS \bf 437\rm. Providence,
R.I., USA, (1990) iv+138 pp.

\bibitem{HJL}
Huntley, J., Jorgenson, J., and Lundelius, R.: On the asymptotic
behavior of counting functions associated to degenerating
hyperbolic Riemann surfaces. \it J. Func. Analysis \bf 149 \rm
(1997) 58-82.

\bibitem{Iwaniec}
Iwaniec, H.: Spectral methods of automorphic forms. Second
edition. Graduate Studies in Mathematics, 53. American
Mathematical Society, Providence, RI; Revista Matemática
Iberoamericana, Madrid (2002) xii+220 pp.

\bibitem{JiZ}
Ji, L. and Zworski, M.: The remainder estimate in spectral
accumulation for degenerating hyperbolic surfaces. \it J. Funct.
Analysis \bf 114 \rm  (1993) 412--420.

\bibitem{Jorgen}
Jorgenson, J.: Asymptotic behavior of Faltings's delta function.
\it Duke Math. J.\bf  61 \rm (1990) 221--254.

\bibitem{JK4}
Jorgenson, J. and Kramer, J.: Bounds for special values of Selberg
zeta functions of Riemann surfaces. \it J. Reine Angew. Math. \bf
541 \rm (2001) 1--28.

\bibitem{JK10}
Jorgenson, J. and Kramer, J.:  Non-completeness of the
Arakelov-induced metric on moduli space of curves.  \it
Manuscripta Math. \bf 119 \rm (2006)

\bibitem{JK11}
Jorgenson, J. and Kramer, J.: Canonical metrics, hyperbolic
metrics, and Eisenstein series for ${{\rm PSL}}_{2}({\bf R})$, in
preparation (unfinished manuscript, (2003)).

\bibitem{JLu1}
Jorgenson, J. and Lundelius, R.: Convergence of the heat kernel
and the resolvent kernel on degenerating hyperbolic Riemann
surfaces of finite volume. \it Quaestiones Math. \bf 18 \rm (1995)
345--363.

\bibitem{JLu2}
Jorgenson, J. and Lundelius, R.: Convergence theorems for relative
spectral functions on hyperbolic Riemann surfaces of finite
volume. \it Duke Math. J. \bf 80 \rm (1995) 785--819.

\bibitem{JoLu}
Jorgenson, J. and Lundelius, R.: Convergence of the normalized
spectral counting function on degenerating hyperbolic Riemann
surfaces of finite volume. \it J. Funct. Anal. \bf 149 \rm (1997)
25--57.

\bibitem{JoLu2}
Jorgenson, J. and Lundelius, R.: A regularized heat trace for
hyperbolic Riemann surfaces of finite volume. \it Comment. Math.
Helv. \bf 72 \rm (1997) 636--659.

\bibitem{JuI}
Judge, C.: Tracking eigenvalues to the frontier of moduli space I:
Convergence and spectral accumulation. \it J. Funct. Analysis \bf
184 \rm (2001) 273--290.

\bibitem{JuII}
Judge, C.:  Tracking eigenvalues to the frontier of moduli space
II: Limits for eigenvalue branches. \it Geom. Funct. Analysis \bf
12 \rm (2002) 93--120.


\bibitem{Kub}
Kubota, T.: Elementary theory of Eisenstein series. Kodansha Ltd.,
Tokyo; Halsted Press [John Wiley \& Sons], New York-London-Sydney,
1973. xi+110 pp.

\bibitem{KM}
Kudla, S. and Millson, J.: Harmonic differentials and closed
geodesics on a Riemann surface. \it Invent. Math. \bf 54 \rm
(1979) 193--211.

\bibitem{Lund}
Lundelius, R.: Asymptotics of the determinant of the Laplacian on
hyperbolic surfaces of finite volume. \it Duke Math. J. \bf 71 \rm
(1993) 211--242.

\bibitem{Pippich}
v. Pippich, A.-M.:  Elliptische Eisensteinreihen, Diplomarbeit,
Humboldt-Universit\"at zu Berlin, (2005).

\bibitem{Randol}
Randol, B.: Cylinders in Riemann surfaces. \it Comment. Math.
Helv. \bf 54 \rm (1979) 1--5.
%\bibitem{WZ}
%Wheeden, Richard L.; Zygmund, Antoni;
%Measure and integral. An introduction to real analysis.
%Pure and Applied Mathematics, Vol. 43. Marcel Dekker, Inc.,
%New York-Basel, 1977. x+274 pp.


\bibitem{Risager}
Risager, M.:  On the distribution of modular symbols for compact
surfaces. \it  International Mathematics Research Notices \bf 41
\rm (2004) 2125--2146.

\bibitem{Wolpert}
Wolpert, S.: Disappearance of cusp forms in special families. \it
Annals of Math. \bf 139 \rm (1994) 239--291.

%\bibitem{WZ} R.~Wheeden, A.~Zymund, \emph{Measure and Integral: An
%Introduction to Real Analysis}, Monographs and texbooks in pure and applied
%mathematics \textbf{43}, Marcel Dekker, Inc., New York and Basel, (1977).
\end{thebibliography}
\end{document}